\documentclass[12pt]{article}

\usepackage{float}
\usepackage{amsmath,amsfonts,amsthm,amssymb,mathtools,bbm} 
\usepackage[colorlinks,citecolor=blue,linkcolor=blue]{hyperref}	
\usepackage{tikz,multirow,multicol,etoolbox} 
\usepackage[inline]{enumitem} 
\usepackage{setspace}
\usepackage{doi}
\usepackage{calc}

\usepackage[longnamesfirst]{natbib}
\bibpunct{(}{)}{;}{a}{,}{,}
\renewcommand{\cite}{\citet}

\DeclareMathOperator{\Z}{\mathbf{Z}}

\DeclareMathOperator{\C}{\mathbf{C}}
\DeclareMathOperator{\D}{\mathbf{D}}
\DeclareMathOperator{\dd}{\mathbf{d}}

\DeclareMathOperator{\X}{\mathbf{X}}
\DeclareMathOperator{\R}{\mathbf{R}}
\DeclareMathOperator{\G}{\mathbf{G}}
\DeclareMathOperator{\y}{\mathbf{y}}

\DeclareMathOperator{\vv}{\mathbf{v}}

\DeclareMathOperator{\I}{\mathbf{I}}
\DeclareMathOperator{\e}{\mathbf{e}}
\DeclareMathOperator{\V}{\mathbf{V}}
\DeclareMathOperator{\Ca}{\mathcal{C}_{\alpha}}
\DeclareMathOperator{\Ma}{\mathcal{M}_{\alpha}}
\DeclareMathOperator{\SIG}{\boldsymbol{\Sigma}}
\DeclareMathOperator{\MU}{\boldsymbol{\mu}}

\DeclareMathOperator{\bta}{\boldsymbol{\beta}}

\DeclareMathOperator{\dlta}{\boldsymbol{\delta}}

\usepackage[utf8]{inputenc}
\usepackage[T1]{fontenc}
\theoremstyle{plain}
\newtheorem{thm}{Theorem}
\newtheorem{lemma}{Lemma}
\newtheorem{prop}{Proposition}

\theoremstyle{definition}

\makeatletter \let\@fnsymbol\@arabic \makeatother 
\usepackage{accents}
\newcommand{\dbtilde}[1]{\accentset{\approx}{#1}}

\usepackage{rotating}

\makeatother
\makeatletter
\renewcommand{\@fnsymbol}[1]{\@arabic{#1}}
\makeatother

\title{Marginal and Conditional Multiple Inference\\
for Linear Mixed Model Predictors}
\author{Peter Kramlinger\thanks{
          peter.kramlinger@univie.ac.at,
          Department of Statistics and Operations Research,
					Universit\"at Wien,
					Oskar-Morgenstern-Platz 1, 1090 Wien, Austria}
\and Tatyana Krivobokova\thanks{
          tatyana.krivobokova@univie.ac.at,
          Department of Statistics and Operations Research,
					Universit\"at Wien,
					Oskar-Morgenstern-Platz 1, 1090 Wien, Austria}
\and Stefan Sperlich\thanks{
          stefan.sperlich@unige.ch,
          Geneva School of Economics and Management,
          Universit\'e de Gen\`eve,
          40 Bd du Pont d'Arve, 1211 Gen\`eve 4, Switzerland}
}

\begin{document}

\baselineskip=25pt

\maketitle

\begin{abstract}

\baselineskip=15pt \noindent
In spite of its high practical relevance, cluster specific multiple inference for linear mixed model predictors has hardly been addressed so far.
While marginal inference for population parameters is well understood, conditional inference for the cluster specific predictors is more intricate.
This work introduces a general framework for multiple inference in linear mixed models
for cluster specific predictors.
Consistent confidence sets for multiple inference are constructed under both, the marginal and the conditional law.
Furthermore, it is shown that, remarkably, corresponding multiple marginal confidence sets are also asymptotically valid for conditional inference.
Those lend themselves for testing linear hypotheses using standard quantiles without the need of re-sampling techniques.
All findings are validated in simulations and illustrated along a study on Covid-19 mortality in US state prisons.
\\
{\textit{Keywords and phrases.}}
Simultaneous inference, multiple testing, mixed parameters, linear mixed models, small area estimation.
\end{abstract}


\baselineskip=20pt
\section{Introduction}
\label{Introduction}
\doublespacing
\linespread{1.2}
\openup -0.22em

Linear mixed models (LMMs) were introduced by 
Henderson in 1950s \citep{Henderson1950, Henderson1953} and are applied if repeated measurements on several independent clusters of interest are available.
The classical LMM, allowing for random intercepts and slopes, is
\begin{equation}
\begin{gathered}\label{LMM}
\y_{i} = \X_{i}\bta + \Z_i \mathbf{v}_i + \e_{i},\;\;i=1,\ldots,m \\
\mathbf{e}_i \sim \mathcal{N}_{n_i}\{\boldsymbol{0}_{n_i},\R_i(\dlta)\}, \quad \mathbf{v}_i \sim \mathcal{N}_{q}\{\boldsymbol{0}_{q},\G(\dlta)\},
\end{gathered}
\end{equation}
with observations $\y_i\in\mathbb{R}^{n_i}$, known $\X_i\in\mathbb{R}^{n_i\times p}$ with $\text{rank}\{(\X_1^t, \dots, \X_m^t)^t\}=p$ for $p\in\mathbb{N}$ fixed and $\Z_i \in\mathbb{R}^{n_i\times q}$, $q\in\mathbb{N}$
fixed, independent random effects $\mathbf{v}_i\in\mathbb{R}^{q}$, and error terms $\e_i\in\mathbb{R}^{n_i}$, such that $\mbox{Cov}(\e_i,\mathbf{v}_i)=\boldsymbol{0}_{n_i\times q}$.
Parameters $\bta\in\mathbb{R}^p$ and $\dlta\in\mathbb{R}^r$, $r\in\mathbb{N}$ fixed, are unknown and we denote $\V_i(\dlta)=\mbox{Cov}(\y_i)=\R_i(\dlta)+\Z_i\G(\dlta)\Z_i^t$, where $\R_i(\dlta)$ and $\G(\dlta)$ are known up to $\dlta$.
LMM (\ref{LMM}) includes the popular nested-error regression model, random coefficient model, and the so-called Fay-Herriot model \citep{Prasad1990}.

Today, LMMs are widely applied in various sciences \citep{tuerlinckx2006statistical, Jiang2007}.
Clusters $i=1, \dots, m$ 
refer for instance to subjects or groups like in biometrics with longitudinal data \citep{LairdEtal1982,LiangZeger1986,verbekeMolenberghs2000},
to treatment levels in medicine \citep{FrancqEtal2019}, or to areas like in the field
of {small area estimation} (SAE), to mention only some prominent application domains.
For the latter, see \cite{TzavidisEtal2018} for a recent review, and \cite{Pratesi2008
}  for examples with interesting time-spatio modeling of $\G(\dlta)$ and $\R_i(\dlta)$.

Depending on the research question, the inference focus may either be on the population parameter $\bta$ or on cluster specific characteristics, and thereby associated with random effects $\mathbf{v}_i$. In the former case, a LMM (\ref{LMM}) can simply be interpreted as a linear regression model with mean $\X_i\bta$ and covariance matrix $\V_i(\dlta)$ that accounts for complex dependence in the data. Inference about $\bta$, where the $\vv_i$ in (\ref{LMM}) are treated as random, is called {\emph{marginal}} and well understood.

Often the interest lies rather in studying {\emph{mixed parameters}}, that is, linear combinations of $\bta$ and $\mathbf{v}_i$, such as $\mu_i = \mathbf{l}_i^t\bta + \mathbf{h}_i^t\mathbf{v}_i$, $i=1,\ldots,m$ with known  $\mathbf{l}_i\in\mathbb{R}^p$ and $\mathbf{h}_i\in\mathbb{R}^{q}$.
In many situations, cf.\ Section 4.1 of \citet{TzavidisEtal2018},
inference about $\mu_i$ with some realized random effects $\mathbf{v}_i$
should then be done {\it conditional} on those $\mathbf{v}_i$,
i.e., $\vv_i$ are treated as fixed.
The importance of this distinction, i.e., between marginal and conditional inference in LMMs, was already emphasized by \cite{Harville1977}, and has attracted particular attention in model selection.
Specifically, \cite{Vaida2005}, who noted that the conventionally used marginal Akaike information criterion is applicable to the selection of $\bta$ only, and suggested a conditional version
for cluster-specific parameters.
Recent contributions \citep{You2016, Lombardia2017} have adopted this distinction and provided bootstrap procedures to accurately estimate the degrees of freedom in the conditional setting.

%

The focus on conditional inference is particularly meaningful if the cluster effects $\mathbf{v}_i$ are rather seen as fixed in practice, and for which
random effects are just a modeling device.
Today, such interpretation is pretty common \citep{Hodges2013}.
Even though it seems then more natural to employ fixed effects models,  in many practical situations their estimators are inefficient, see e.g. \citet{pfeffermann2013}.
Then, a reasonable approach to obtain estimators of cluster specific effects is to employ model (\ref{LMM}) as if $\mathbf{v}_i$ were random and obtain a predictor for them.
Yet, to perform inference as if $\mathbf{v}_i$ were fixed, one needs to condition on the cluster, i.e., on $\mathbf{v}_i$. 
For example, in the application in Section \ref{RealStud} the mortality in US state prisons is studied. The effects $\mathbf{v}_i$ in this example model the state specific effects on Covid-19 mortality (e.g., due to state policy and/or population structure).
Since too few observations per state are available, the state effect is predicted within the modeling framework of a LMM.
Assume that one is interested in inference about the mean mortality in each state, that is, including the state effect $\mathbf{v}_i$. Since the state effect is not necessarily considered to be random in nature and the inference focus is on the state level, the corresponding  inference should be conducted under the conditional law, that is, treating $\mathbf{v}_i$ as fixed.
Therefore, the main focus of this work is on conditional inference on $\mu_i$. For more discussion on marginal versus conditional inference see also our Supplement, Section 6.1. 


There is a large body of literature on constructing confidence intervals for each $\mu_i$ separately under the marginal law. Since under the marginal law, estimators $\hat\mu_i$ of $\mu_i$ obtained from (\ref{LMM}) are unbiased, much attention has been given to the estimation of the mean squared error $\mbox{MSE}(\hat{\mu}_i) = \mbox{E}(\mu_i-\hat{\mu}_i)^2=\mbox{Var}(\hat\mu_i)$, where the expectation is taken under the marginal law, that is, treating $\vv_i$ as random. To estimate this marginal MSE, one can either plug in an appropriate estimator of $\dlta$,
or use unbiased marginal MSE approximations \citep{Prasad1990, Datta2000, Das2004}. Other approaches to the estimation of marginal MSE comprise a diverse collection of bootstrap methods \citep{Manteiga2008, Chatterjee2008}.
%
The conditional inference on single $\mu_i$, that is conditioning on $\vv_i$, turns out to be infeasible due to the bias of $\hat\mu_i$ which arises under the conditional law. While estimation of the bias leads to unacceptably wide intervals \citep{Datta2002, Jiang2006}, ignoring it leads to strong under-coverage. This was mentioned as an open problem by \cite{pfeffermann2013}.

Conditional and marginal inference about all $\mu_1,\ldots,\mu_m$ simultaneously or about a subset thereof  has been largely neglected. To the best of our knowledge, only \cite{Ganesh2009} considered a related problem of Bayesian inference about certain linear combinations of $\mu_i$ in the Fay-Herriot model.
\citet{Kasia2019} and \cite{Kasia2019b} used max-type statistics to construct simultaneous intervals for mixed parameters $\mu_i$ of generalized LMM under the marginal law. 
For a discussion of the average coverage of cluster specific confidence intervals see \cite{Zhang2007}, and Section \ref{SimStud} for its relation to our method.
None of these contributions considered multiple inference under the conditional law.

%

Altogether, there is a lack of results on multiple inference in linear mixed models and a tension between marginal and conditional focus in inference. In this work we address both issues.
First, we construct confidence sets for $\mu_1,\ldots,\mu_m$ in LMMs. Second, we consider those joint (or multiple) confidence sets under both, the conditional and the marginal law.
For the former we show that the nominal coverage is attained at the usual parametric rate. Then we show that, surprisingly, joint  confidence sets built under the marginal law, are accurate at the same parametric rate, and also approximately valid when conditioning on the clusters. This, however, is not true in general for the cluster-wise confidence intervals, i.e.,  for single $\mu_i$.
Next, we use the derived confidence sets to develop multiple tests for linear hypotheses, both on all $\mu_1, \dots, \mu_m$ or on a subvector thereof.
Finally, the practical use and relevance of the derived methods is illustrated in simulations and a study on Covid-19 mortality in US state prisons.

The main results are given in Section \ref{SimInf} including applications for comparative statistics and testing linear hypotheses. The results are visualized via simulations in Section \ref{SimStud} and a practical application is given in Section \ref{RealStud}. We conclude with a discussion in Section \ref{Discussion}. Relevant proofs are deferred to the Appendix, while auxiliary proofs and additional results are provided in the Supplement.

\section{Confidence Stets for Multiple Inference}\label{SimInf}

\subsection*{Marginal Simultaneous Prediction Sets}

We start by introducing further notation and assumptions, for a general
monograph on LMMs and generalizations see e.g. 
\cite{Demidenko2004}.
For model (\ref{LMM}), under the marginal law, the best linear unbiased predictor (BLUP) of $\mu_i = \mathbf{l}_i^t\bta + \mathbf{h}_i^t\mathbf{v}_i$ reads as
\begin{equation}
\begin{gathered}\label{BLUP}
\tilde{\mu}_i = \tilde{\mu}_i\left\{\dlta, \hat\bta\big(\dlta\big) \right\}
= \mathbf{l}_i^t\hat\bta\big(\dlta\big) + \mathbf{b}_i(\dlta)^t\left\{\y_i-\X_i\hat\bta\big(\dlta\big)\right\};\vspace*{5pt}\\ \mbox{where } \
\mathbf{b}_i(\dlta)^t = \mathbf{h}_i^t\G(\dlta)\Z_i^t{\V}_i(\dlta)^{-1},
\ \mbox{ and} \\
\hat\bta\big(\dlta\big)  = \left\{\sum_{i=1}^m \X_i^t{\V}_i(\dlta)^{-1}\X_i\right\}^{-1}\sum_{i=1}^m \X_i^t{\V}_i(\dlta)^{-1}\y_i.
\end{gathered}
\end{equation}
%
If the variance components $\dlta$ are unknown, they can be estimated using restricted maximum likelihood (REML) 
as given in (15) in the Supplement, or by Henderson III,
as defined by \cite[Chapter 5]{Searle1992}. 
%
Replacing $\dlta$ in (\ref{BLUP}) by an estimator based on either one of these methods, gives the empirical BLUP (EBLUP)
\begin{align} \label{EBLUP}
\hat{\mu}_i = \tilde{\mu}_i\left\{\hat\dlta, \hat\bta\big(\hat\dlta\big) \right\}.
\end{align}

Subsequently, the dependency on $\hat\dlta$ or $\dlta$ is suppressed if it is clear from the context.
%
%
Consider the asymptotic scenario 
\begin{enumerate}[label=(A\arabic*), itemsep=2pt]    \vspace{-0.5cm}
\item  $m\rightarrow \infty$ while $\sup_i n_i = O(1)$. \hspace{1cm} \label{SAE}
\end{enumerate} \vspace{-0.5cm}
It encompasses the standard SAE assumption: For a growing number of clusters there are few observations per cluster.
The requirement that $m\rightarrow \infty$ ensures consistent estimation of $\hat\bta(\dlta)$ \citep[Section 3.6.2]{Demidenko2004}.
The boundedness condition on the cluster sample sizes is not crucial for the results we derive subsequently, but rather constitute the most unfavorable case under which they hold true.
In particular, if some or all $n_i\rightarrow\infty$, certain rates may only improve, for more details see the discussion in the Appendix.

%
Further, we work with the quite standard, though adapted, regularity conditions
\begin{enumerate}[label=(B\arabic*), itemsep=2pt]  \vspace{-0.5cm}
\item $\X_i$, $\Z_i$, $\G(\dlta)>0$, $\R_i(\dlta)>0$, $i=1,\ldots,m$ contain only bounded values. \label{A}
\item $\mathbf{d}_i^t = \mathbf{l}_i^t - \mathbf{b}_i(\dlta)^t\X_i$ has entries $d_{ik} = O(1)$ for $k= 1, \dots, p$.\label{C}
\item $\big\{\frac{\partial}{\partial\delta_j}\mathbf{b}_i(\dlta)^t\X_i\big\}_k = O(1)$, for $j=1,\dots,r$ and $k = 1, \dots, p$.\label{D}
\item $\V_i(\dlta)$ is linear in the variance components $\dlta$.\label{E}
\end{enumerate}  \vspace{-0.5cm}

The last condition \ref{E} implies that the second derivatives of $\R_i$ and $\G$ w.r.t. $\dlta$ are zero.
These assumptions imply that $\mbox{E}(\hat{\mu}_i-\mu_i)=0$  \citep{Jiang2000}.


Subsequently, dropping the cluster index $i$ refers to the respective quantity over all clusters: $\y = (\y_1^t,\dots,\y_m^t)^t$, $\V(\dlta) = \text{diag}\{\V_i(\dlta)\}_{i=1 ,\dots,m}$, $\X = (\X_1^t, \dots, \X_m^t)^t$, etc.
Now we can construct prediction sets for $\MU=(\mu_1,\ldots,\mu_m)^t$.
The theory below is based on an extension of the MSE estimator for point-wise marginal inference from \cite{Prasad1990} for multiple inference.
We start by constructing a prediction set $\Ma$ such that $\mbox{P}(\MU\in \Ma) \approx 1-\alpha$, where $\mbox{P}$ refers to the marginal probability under (\ref{LMM}), for a pre-specified level $\alpha\in(0,1)$. That is, $\MU$ is considered as a random variable under the marginal law and the corresponding prediction sets are not meant for the conditional inference about a fixed $\MU$.
Consider an estimator for $\SIG = \mbox{Cov}(\hat{\MU}-\MU)$ given by
\begin{eqnarray} \label{SIGorig}
\SIG &=& \mathbf{K}_1({\dlta}) + \mathbf{K}_2({\dlta}) + \mathbf{K}_3({\dlta})
\ , \ \mbox{ with  }
\\ \nonumber
\mathbf{K}_1({\dlta})
&=& \mbox{Cov}\left(\dbtilde{\MU}-\MU\right)
= \text{diag}\left[ \mathbf{h}_i^t\left\{\G-\G\Z_i^t\V_i^{-1}\Z_i\G\right\}\mathbf{h}_i   \right]_{i=1,\dots,m},\\ \nonumber
\mathbf{K}_2({\dlta})
&=& \mbox{Cov}\left(\tilde{\MU}-\dbtilde{\MU}\right)
= \left\{\mathbf{d}_i^t\left(\sum_{l=1}^m \X_l^t\V_l^{-1}\X_l\right)^{-1}\mathbf{d}_k\right\}_{i,k=1,\dots,m},
\\ \nonumber
\mathbf{K}_3({\dlta})
&=& \mbox{Cov}\left(\hat{\MU}-\tilde{\MU}\right),
\end{eqnarray}
where
$\mathbf{d}_i$ as in  \ref{C} and
$\dbtilde{\MU} = (\dbtilde{\mu}_1, \dots, \dbtilde{\mu}_m)^t$ with $\dbtilde{\mu}_i = \tilde{\mu}_i(\dlta, \bta) 
$.
The decomposition (\ref{SIGorig}) partly follows results of \cite{Kackar1984}. 
The following lemma gives an estimator for $\SIG$, and evaluates its bias, which will be needed later on.
\begin{lemma} \label{marglemma}
Let model (\ref{LMM}) hold and $\hat{\dlta}$ be either a REML estimator or given by Henderson III. Under \ref{SAE} and \ref{A}-\ref{E}, consider the estimator $\widehat{\SIG} = \widehat\SIG(\hat\dlta)$ for $\SIG$ given by
\begin{equation}
\begin{aligned} \label{defSIGMA}
\widehat\SIG(\hat\dlta) &= \mathbf{K}_1({\hat\dlta}) + \mathbf{K}_2({\hat\dlta}) + 2\widehat{\mathbf{K}}_3({\hat\dlta}),\\
\widehat{\mathbf{K}}_3({\hat\dlta})
&= \text{diag}\left[\text{tr}\left\{\frac{\partial\mathbf{b}_i^t}{\partial\dlta}
\V_i\frac{\partial\mathbf{b}_i}{\partial\dlta^t}\overline{\V}\right\}\right]_{i=1,\dots,m},
\end{aligned}
\end{equation}
	where $\overline{\V}$ is the asymptotic covariance matrix of $\hat\dlta$. It then holds
\begin{align*}
\mbox{E}\big(\widehat\SIG\big) = \SIG + \big\{O\big(m^{-3/2}\big)\big\}_{m\times m}.
\end{align*}
\end{lemma}
Here, $\{O(m^{-3/2})\}_{m\times m}$ denotes an $(m \times m)$ matrix with each entry being of order $O(m^{-3/2})$.
This error term comes from the uncertainty in estimating $\dlta$.
The result can be concluded from the point-wise case, which corresponds to the diagonal entries of $\widehat\SIG$, and was shown by \cite{Prasad1990} for Henderson III and by \cite{Datta2000} for REML.
Since only $\mathbf{K}_2(\hat\dlta)$ contributes to off-diagonal entries, Lemma \ref{marglemma} is a straightforward, though a tedious extension, so that we skip its proof.
With $\widehat\SIG$ at hand we can state now:
\begin{thm} \label{margdist}
Let model (\ref{LMM}) hold and $\widehat{\SIG} = \widehat\SIG(\hat\dlta)$ as given in (\ref{defSIGMA}). Under \ref{SAE} with \ref{A}-\ref{E} it holds that
\begin{align*}
\mbox{P}\bigg\{\big\|\widehat{\SIG}^{-1/2}(\hat{\MU}-\MU)\big\|^2 < \chi^2_{m,1-\alpha}\bigg\} = 1-\alpha  + O(m^{-1/2}),
\end{align*}
where $\alpha\in(0,1)$ and $\chi^2_{m,1-\alpha}$ is the $\alpha$-quantile of the $\chi_m^2$-distribution.
\end{thm}
%
As before, the error rate is due to the uncertainty in estimating $\dlta$.
The error rates with Lemma \ref{marglemma} and Theorem \ref{margdist} differ, since the theorem simultaneously considers all $m$ elements of $\MU\in\mathbb{R}^m$, so that the error is increased from $O(m^{-3/2})$ to $O(m^{-1/2})$.
From Theorem \ref{margdist} we immediately obtain the prediction set under the marginal law,
\begin{align*}
\Ma = \bigg\{ \MU\in\mathbb{R}^m: \big\|\widehat{\SIG}^{-1/2}(\hat{\MU}-\MU)\big\|^2 \leq \chi^2_{m,1-\alpha} \bigg\},
\end{align*}
with $P(\MU\in \Ma)\approx 1-\alpha$, for $\alpha\in(0,1)$.
In the marginal case, $\Ma$ is a prediction region for the random variable $\MU$ and can therefore not be readily interpreted as a confidence region for a fixed $\MU$.

\subsection*{Conditional Simultaneous Confidence Sets}

If the inference focus is conditional, i.e.,\ $\mathbf{v}$ is treated as fixed,
then the aim is to obtain a confidence set $\Ca$ with  $\mbox{P}(\MU\in \Ca | \mathbf{v}) \approx 1-\alpha$.
Notation $\text{P}(\cdot|\mathbf{v})$ means that the probability is taken under model (\ref{LMM}), in which $\mathbf{v}= (\mathbf{v}^t_1, \dots, \mathbf{v}_m^t)^t$ is given,
or `fixed'.
Since small cluster sample sizes result in unreliable direct estimators, the confidence set $\Ca$ is based on the EBLUP $\hat\mu_i$ for $\mu_i$ from (\ref{EBLUP}).

If effect $\mathbf{v}$ is some fixed parameter, not necessary a realization of a random variable, then model (\ref{LMM}) is misspecified and $\dlta$ is not meaningful. Since the parameters are still estimated from (\ref{LMM}), one needs to replace $\dlta$ by $\dlta^v$ which is an oracle parameter that we define by $\dlta^v = \mbox{E}(\hat\dlta|\mathbf{v})$.
To control for its variation, $\mathbf{v}$ needs to meet the following conditions which are rather general:
\begin{enumerate}[label=(C\arabic*), itemsep=2pt]   \vspace{-0.5cm}
\item $\sum_{i=1}^m (\vv_i)_{e} = O(m^{1/2})$, $e = 1,\dots, q$; \label{C1}
\item $\sum_{i=1}^m \{\vv_i\vv_i^t - \G(\dlta^v)\}_{ef} = O(m^{1/2})$, $e,f = 1,\dots, q$. \label{C2}
\end{enumerate}   \vspace{-0.5cm}
The first condition is required to identify $\bta$ from $(\bta^t, \vv^t)$.
Variants thereof are commonly used in econometrics, see \citet[Section 3.2]{Hsiao2014}.
Condition \ref{C2} states that the $\mathbf{v}_i$'s should not be too different from each other; in particular, it ensures that the stochastic part of the observed information matrix
is dominated by its deterministic part.
A formal discussion and details are given in Lemma 4 in the Supplement.
If $\mathbf{v}$ is a realization of a normally distributed random variable, then conditions \ref{C1} and \ref{C2} are readily satisfied.

As in the marginal case, consider the standard regularity conditions \ref{A}-\ref{D} and asymptotic scenario \ref{SAE}.
The latter implies that under conditional law $\mbox{E}(\hat\MU - \MU|\mathbf{v}) \nrightarrow \mathbf{0}_m$, due to the boundedness of  $n_i$, rendering conditional inference for single $\mu_i$ infeasible.
If $m\rightarrow\infty$ and $n_i\rightarrow\infty$ for some fixed $i$, such inference would be possible for corresponding $\mu_i$, as $\mbox{E}(\hat\mu_i - \mu_i |\mathbf{v})\rightarrow 0$.
Only if $m\rightarrow\infty$ and $n_i\rightarrow\infty$ for all $i = 1, \dots, m$, the conditional bias vanishes for all clusters and $\mbox{E}(\hat\MU - \MU|\mathbf{v}) \rightarrow \mathbf{0}_m$.
Since our results allow multiple conditional inference for non-vanishing bias under \ref{SAE}, they still holds also if all or some $n_i \rightarrow\infty$.
These technical differences are discussed in the Appendix in more detail.

Proceeding as for the marginal case, $\SIG_v=\mbox{Cov}(\hat{\MU}-\MU|\mathbf{v})$ can be decomposed as
\begin{align} \label{SIGCorig}
\SIG_v = \mathbf{L}_1({\dlta^v}) + \mathbf{L}_{2}({\dlta^v}) + \mathbf{L}_3({\dlta^v}) + \mathbf{L}_4({\dlta^v});
\end{align}
where for $\mathbf{K}_k = \mathbf{R}_k\V_k^{-1}\X_k\left(\sum_{l=1}^m \X_l^t\V_l^{-1}\X_l\right)^{-1}$, and with notation from (\ref{BLUP})
\begin{eqnarray*}
\mathbf{L}_1({\dlta^v})
&=& \mbox{Cov}\left(\dbtilde{\MU}\big| \mathbf{v}\right)
= \text{diag}\left\{ \mathbf{b}_i^t\R_i\mathbf{b}_i   \right\}_{i=1,\dots,m},\\
\mathbf{L}_2({\dlta^v})
&=& \mbox{Cov}\left(\tilde\MU\big| \mathbf{v}\right)  - \mbox{Cov}\left(\dbtilde\MU\big| \mathbf{v}\right)
= \left\{
\mathbf{b}_k^t\mathbf{K}_k\mathbf{d}_i
+\mathbf{b}_i^t\mathbf{K}_i\mathbf{d}_k
+ \sum_{l=1}^m\mathbf{d}_i^t \mathbf{K}_l^t\mathbf{R}_l^{-1}\mathbf{K}_l\mathbf{d}_k
\right\}_{i,k=1,\dots,m}\\
\mathbf{L}_3({\dlta^v})
&=& \mbox{Cov}\left(
\hat\MU-\tilde\MU, \tilde\MU
\big| \mathbf{v}\right) + \mbox{Cov}\left(
\hat\MU-\tilde\MU, \tilde\MU
\big| \mathbf{v}\right)^t,\\
\mathbf{L}_4({\dlta^v})
&=& \mbox{Cov}\left(  \hat\MU-\tilde\MU \big| \mathbf{v}\right)  .
\end{eqnarray*}
This decomposition is similar to (\ref{SIGorig}), but the cross-terms do not vanish.
As explained above, $\dlta^v$ substitutes now $\dlta$, although $\hat{\dlta}$ remains the same.
The next lemma gives an estimator $\widehat\SIG_v= \widehat\SIG_v(\hat\dlta)$ for the conditional covariance matrix $\SIG_v$, and evaluates its bias.
\begin{lemma}\label{LemmaCondReml}
Let model (\ref{LMM}) hold. Under \ref{SAE}, with \ref{A}-\ref{E}, \ref{C1}, \ref{C2},
\begin{align} \label{defcSIGMA}
\widehat\SIG_v(\hat\dlta) &=
\mathbf{L}_1(\hat{\dlta})
+ \mathbf{L}_2(\hat{\dlta})
+ \widehat{\mathbf{L}}_3(\hat{\dlta})
+ \widehat{\mathbf{L}}_4(\hat{\dlta})
- \widehat{\mathbf{L}}_5(\hat{\dlta}),
\end{align}
where $\widehat{\mathbf{L}}_3(\dlta^v)$ is given in (\ref{L3REML})
if $\hat{\dlta}$ is obtained via REML,  or in (\ref{L3H3}) if $\hat{\dlta}$ is obtained via Henderson III. Further,
\begin{align*}
\widehat{\mathbf{L}}_4({\dlta^v})
&= \text{diag}\bigg[\text{tr}\bigg\{\frac{\partial\mathbf{b}_i^t}{\partial\dlta^v}
\R_i\frac{\partial\mathbf{b}_i}{\partial(\dlta^v)^t}\overline{\V}\bigg\}
\bigg]_{i=1,\dots,m}, \\
\widehat{\mathbf{L}}_5({\dlta^v})
&= \frac{1}{2}\text{diag}\bigg\{\text{tr}\bigg[ \frac{\partial^2\{\mathbf{L}_1(\dlta^v)\}_{ii}}{\partial\dlta^v\partial(\dlta^v)^t}  \overline{\V}\bigg]\bigg\}_{i=1,\dots,m},
\end{align*}
where $\overline{\V}$ is the asymptotic covariance matrix of $\hat\dlta$.
Then it holds
\begin{align*}
\mbox{E}\big(\widehat\SIG_v\big|\mathbf{v}\big) = \SIG_v + \big\{O\big(m^{-3/2}\big)\big\}_{m\times m}.
\end{align*}
\end{lemma}
The proof is given in the Supplement.
As the EBLUP is not unbiased under conditional law, $(\hat\MU - \MU)^t\SIG_v^{-1}(\hat\MU - \MU)|\mathbf{v}\sim \chi^2_m(\lambda)$ for $\lambda = \| \SIG_{v}^{-1/2} \mbox{E}(\hat\MU - \MU|\mathbf{v}) \|^2$.
The non-centrality parameter $\lambda$ depends on the conditional bias, and cannot be estimated  directly for any cluster individually, but only jointly. Specifically, let
${\mathbf{A}}(\dlta^v) = ({\mathbf{a}}_1^t,\dots, {\mathbf{a}}_m^t)^t\in\mathbb{R}^{m\times n}$, with
\begin{equation} \label{defA}
{\mathbf{a}}_i = (\mathbf{b}_i^t\mathbf{Z}_i - \mathbf{h}_i^t)\mathbf{J}_i(\Z^t\Z)^{-1}\Z^t + \mathbf{d}_i^t(\X^t\V^{-1}\X)^{-1}\X^t\V^{-1},
\end{equation}
where $\mathbf{J}_i=(0,\dots, 0, \mathbf{I}_{q}, 0,\dots, 0)\in\mathbf{R}^{q\times qm}$, so that ${\mathbf{a}}_i\mathbf{Z}\mathbf{v} = \mbox{E}(\tilde\mu_i-\mu_i|\mathbf{v})$.
We propose to estimate the non-centrality parameter by
\begin{eqnarray} \label{lambdahat}
& \hat\lambda = \max\left[ 0, \tilde\lambda\big\{ \widehat\SIG_v(\hat\dlta),\hat\bta, \hat\dlta \big\} \right],  & \\ &  \nonumber
\tilde\lambda\big( \SIG_v, \bta, \dlta^v \big)
= \big\|\SIG_v^{-1/2}\mathbf{A}(\dlta^v)\y\big\|^2 - \big\|\SIG_v^{-1/2}\mathbf{A}(\dlta^v)\mathbf{R}(\dlta^v)^{1/2}\big\|^2 - \big\|\SIG_v^{-1/2}\mathbf{A}(\dlta^v)\mathbf{X}\bta\big\|^2.
\end{eqnarray}
Note that $\text{E}\{ \tilde\lambda\big( \SIG_v, \tilde\bta, \dlta^v \big) |\mathbf{v}  \}/n =  \lambda + O(m^{-1/2})$.
With this estimator we can show
\begin{thm} \label{conddist}
Let model (\ref{LMM}) hold and $\widehat{\SIG}_v = \widehat{\SIG}_v(\hat\dlta)$ as in (\ref{defcSIGMA}) and $\hat\lambda$ as in (\ref{lambdahat}). Under \ref{SAE}, with \ref{A}-\ref{E}, \ref{C1} and \ref{C2} it holds that
\begin{align*}
\mbox{P}\bigg\{\big\|\widehat{\SIG}_v^{-1/2}(\hat{\MU}-\MU)\big\|^2 < \chi^2_{m,1-\alpha}(\hat\lambda)  \bigg|\mathbf{v}\bigg\} = 1-\alpha  + O(m^{-1/2}),
\end{align*}
where $\alpha\in(0,1)$, and $\chi^2_{m,1-\alpha}(\hat\lambda) $ is the 
$\alpha$-quantile of the non-central $\chi_m^2$-distribution.
\end{thm}
Like in Theorem \ref{margdist}, the error of rate $m^{-1/2}$ is due to the uncertainty in estimating $\dlta^v$ which enters $\widehat\SIG_v$ and $\hat\lambda$.
The result gives the conditional confidence set
\begin{align*}
\Ca = \bigg\{ \MU\in\mathbb{R}^m: \big\|\widehat{\SIG}_v^{-1/2}(\hat{\MU}-\MU)\big\|^2 \leq \chi^2_{m,1-\alpha}(\hat{\lambda}) \bigg\},
\end{align*}
with $\mbox{P}(\MU\in \Ca | \mathbf{v}) \approx 1-\alpha$, $\alpha\in(0,1)$. The practical difficulty when constructing $\Ca$ is the unhandy representation of  $\widehat\SIG_v$ and $\hat\lambda$.
Yet, the following result states that the much simpler $\Ma$, albeit derived for the marginal case, leads to the asymptotically correct coverage in the conditional case as well.
\begin{thm} \label{margforcond}
Let model (\ref{LMM}) hold and $\widehat{\SIG} = \widehat\SIG(\hat\dlta)$ as in (\ref{defSIGMA}). Under \ref{SAE} with \ref{A}-\ref{E}, \ref{C1} and \ref{C2} it holds that
\begin{align*}
\mbox{P}\bigg\{\big\|\widehat{\SIG}^{-1/2}(\hat{\MU}-\MU)\big\|^2 < \chi^2_{m,1-\alpha}  \bigg|\mathbf{v}\bigg\} = 1-\alpha  + O(m^{-1/2}) .
\end{align*}\end{thm}
The theorem states that the misspecification in using the marginal covariance matrix under the conditional law is averaged out across clusters.
Remarkably, the rate at which this misspecification vanishes is of the same magnitude as the estimation error in estimating $\dlta^v$.
Without any extra cost, at least if looking at the first order, the much simpler marginal confidence set can be applied in the conditional scenario, i.e., $\mbox{P}(\MU\in \Ma | \mathbf{v}) \approx 1-\alpha$.

In contrast to the previous two theorems, the error term in Theorem \ref{margforcond} is composed by both the estimation error, which relies on all observations, and the misspecification error, which relies on the number of comparisons.
This highlights why individual conditional inference based on the marginal MSE is not possible under \ref{SAE}:
If the quadratic form in Theorem \ref{margforcond} is reformulated for one cluster, it follows from the proof that
\begin{equation} \label{cond-clusterwise}
\mbox{P}\bigg\{\frac{(\hat{\mu}_i-\mu_i)^2}{\hat{\sigma}_{ii}} < \chi^2_{1,1-\alpha}  \bigg|\mathbf{v}\bigg\} = 1-\alpha  + O\left(n_i^{-1/2}\right),
\end{equation}
which is not useful for $\mbox{sup}_i n_i = O(1)$.
Only if $n_i\rightarrow\infty$, $\hat\mu_i$ becomes consistent for $\mu_i$ under the conditional law, and nominal coverage for a single $\mu_i$ is asymptotically attained.

\begin{figure}[t]
	\centering
	\begingroup
	\tikzset{every picture/.style={scale=0.8}}%
  \input{sketch.tex}
	\endgroup
	\caption{
		Sketch comparing $\Ca$ (dashed) and $\Ma$ (solid) as confidence ellipses for $m=2$.
	}
	\label{sketch}
\end{figure}

Theorems \ref{conddist} and \ref{margforcond} deal with the problem on how conditional inference for mixed parameters could be performed.
The latter theorem suggests that multiple inference about $\MU$ under the conditional law can be performed based on the confidence sets obtained under the marginal law.
Figure \ref{sketch} shows that this effect occurs even though the sets are not necessarily equal.
For $m=2$, the two confidence sets are drawn for randomly generated random effects as described in Section \ref{SimStud}.
Although being centered around $\hat\MU$ and holding the same coverage probability under the conditional law, they differ in shape.
Besides the obvious dependence of both $\Ca$ and $\Ma$ on the random effect $\mathbf{v}$ through $\hat\dlta$, the former also depends on $\mathbf{v}$ via $\lambda$.
The non-centrality parameter extends $\Ca$ to such a degree that the confidence region meets nominal coverage probability.
The set $\Ma$, on the other hand, ignores the bias of $\hat\MU-\MU$ that occurs under the conditional law.
But this is compensated in that it is inflated by the marginal variance of $\hat\MU-\MU$, which, in contrast to the conditional variance incorporates the variability of the random effects.
Theorem \ref{margforcond} postulates that both properties cancel each other out, in such a manner that nominal level is approximately attained.
The obvious suspicion that this occurs at the cost of a larger volume for the marginal set
gets dispelled in our simulations in Section \ref{SimStud}.
Clearly, $\Ma$ does not require any separately estimated parameters, which simplifies its implementation.


\subsection*{Conditional Multiple Testing} 

It is appealing to use the derived results for multivariate hypothesis testing under the conditional law.
This can be used to see if $\MU$ lies in a given subspace of $\mathbb{R}^m$, and includes tests of all kind of linear comparisons between clusters. It can also be applied to examine if cluster specific effects are present within subsets, cf.\ Section \ref{RealStud}.
%
Consider
\begin{equation}
\label{eq:test}
H_0: \; \mathbf{L}(\MU - \mathbf{a}) = \boldsymbol{0}_{u} \quad\mbox{vs.}\quad H_1:\quad \mathbf{L}(\MU - \mathbf{a}) \neq \boldsymbol{0}_{u} ,
\end{equation}
where $\mathbf{a}\in\mathbb{R}^m$ and $\mathbf{L}$ is a given $(u \times m)$-matrix with $u\leq m$ and $\text{rank}(\mathbf{L}) = u$, $u = m^{\xi_1}$, and $\xi_1\in(0,1]$ bounded away from zero.
The dimension $u$ of the linear subspace of $\mathbb{R}^m$ corresponds to the number of multiple tests of linear combinations, whereas each linear combination of interest is specified in the rows of $\mathbf{L}$.
For example, for $\mathbf{L} = \mathbf{I}_m$ and $\mathbf{a} = (a_1, \dots, a_m)^t$, $a_i \neq a_j$, $i,j\leq m$, tests whether the mixed parameters take on some ex-ante assumed value(s).
For conditional inference about $\MU$, Theorem \ref{conddist} gives the $\alpha$-level test for (\ref{eq:test}), that rejects $H_0$ if $\mathbf{a} \not\in \mathcal{C}_{\alpha, \mathbf{L}} $, where
\begin{align*}
\mathcal{C}_{\alpha, \mathbf{L}} = \left\{ \mathbf{a}\in\mathbb{R}^{m}: \big\|\big(\mathbf{L}\widehat{\SIG}_v\mathbf{L}^t\big)^{-1/2}\mathbf{L}(\hat\MU - \mathbf{a}) \big\|^2 \leq \chi^2_{u,1-\alpha}(\hat\lambda_{\mathbf{L}})\right\}.
\end{align*}
with $\hat\lambda_{\mathbf{L}}$ being the non-centrality parameter estimate that depends on the covariance $\mathbf{L}\widehat{\SIG}_v\mathbf{L}^t$. This test is consistent with an error of size $O(u^{-1/2})$.
Theorem \ref{margforcond} allows us to employ the confidence set $\Ma$ as well. An $\alpha$-level test rejects $H_0$ if $\mathbf{a} \not\in \mathcal{M}_{\alpha, \mathbf{L}}$, where
$$
\mathcal{M}_{\alpha, \mathbf{L}} = \left\{ \mathbf{a}\in\mathbb{R}^{m}: \big\|\big(\mathbf{L}\widehat{\SIG}\mathbf{L}^t\big)^{-1/2}\mathbf{L}(\hat\MU - \mathbf{a}) \big\|^2 \leq \chi^2_{u,1-\alpha} \right\}.
$$
This test is again consistent at rate $u^{-1/2}$.
It affirms that individual confidence intervals ($u=1$) cannot be constructed using neither $\mathcal{M}_{\alpha, \mathbf{L}}$ nor $\mathcal{C}_{\alpha, \mathbf{L}}$ under \ref{SAE}, the standard SAE assumption.
Note finally that the derived estimators $\widehat\SIG$ and $\widehat\SIG_v$ also lend themselves for related testing procedures, such as Tukey's tests, see our Supplement.

\section{Simulation Examples and Performance Study} \label{SimStud}

While it is to be emphasized that the above methods, theory and developments hold for the general LMM (\ref{LMM}), and thereby also for complex models with slopes that potentially vary over clusters, we concentrate in the following on a most popular though simpler version, namely the
nested error regression model from \cite{Battese1988} with
$e_{ij}\sim\mathcal{N}(0,\sigma_e^2),$ $v_i\sim\mathcal{N}(0,\sigma_v^2)$,
and
\begin{equation}
\begin{gathered}\label{NER}
y_{ij} = \beta_0 + x_{ij}\beta_1 + v_i + e_{ij},\;\;i=1,\ldots,m,\;j=1,\ldots,n_i .
\end{gathered}
\end{equation}
%
The data are simulated as follows.
For each given set of the parameters $m$, $n_i$, $\sigma_e^2$, $\sigma^2_v$, the value of the cluster effect $v_i$ is obtained as a realization of a $\mathcal{N}(0,\sigma^2_v)$ distributed random variable and remains fixed in all Monte Carlo samples.
The covariates $x_{ij}$ are drawn once from a standard normal distribution, whereas the coefficient parameters are set to $(\beta_0, \beta_1)  = (-4.9, 0.03)$, which is similar to the study in Section \ref{RealStud}.
The parameter of interest is the conditional mean $\mu_i = \beta_0 + \sum_{j=1}^{n_i} x_{ij} \beta_1 / n_i + v_i$.
Since the random effects are being drawn from a Gaussian distribution, the requirements of Theorem \ref{margforcond} are fulfilled.

\begin{figure}[t]
 	    	\centering
 		\input{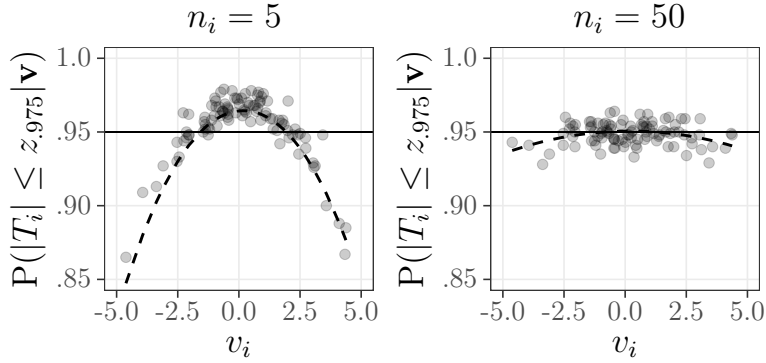}
 		\caption{
		Where inference fails:
		Empirical coverage of marginal 95\% cluster-wise confidence intervals for $\mu_i$ under the conditional law for small (left) and large (right) cluster effects do not meet nominal level. The dashed lines give the theoretical coverage. }
 		\label{cov_fig}
 \end{figure}

Before we study the joint inference, let us briefly look at the cluster-wise one that is, about specific $\mu_i$. As (\ref{cond-clusterwise}) indicates that this cannot be done consistently under that conditional law if $n_i \nrightarrow \infty$, one considers statistics of type
$T_i := (\tilde{\mu}_i-\mu_i)\mbox{Var}(\tilde{\mu}_i-\mu_i)^{-1/2}$, $i=1,...,m$
under the marginal law.
Yet, plotting $\mbox{P}(|T_i|\leq z_{1-\alpha/2}|\mathbf{v})$ for all $v_i$, with $z_{1-\alpha/2}$ being the two-sided $\alpha$-quantile of $\mathcal{N}(0,1)$, Figure \ref{cov_fig} shows how much the coverage probabilities of standard confidence intervals vary with the cluster effects $v_i$.
These results are based on 1.000 Monte Carlo samples with $m=100$, $(\sigma_v^2, \sigma_e^2) = (4,4)$. We see that clusters which comprise a large $|v_i|$, i.e., with most prominent cluster effect, exhibit a severe undercoverage. This is particularly annoying, since such clusters are arguably those that a practitioner might be most interested in \citep{Jiang2006}.  For large $n_i$, this problem is less pronounced, since the bias for every cluster vanishes asymptotically, and so does the difference between conditional and marginal variance.

On average, i.e., over all clusters, over- and undercoverage cancel out each other, as the following result shows.
\begin{prop} \label{avgcov}
	Let model (\ref{LMM}) hold, $\dlta^v$ known, $T_i = (\tilde{\mu}_i-\mu_i)\mbox{Var}(\tilde{\mu}_i-\mu_i)^{-1/2}$ and $z_{1-\alpha/2}$ the two-sided $\alpha$-quantile of $\mathcal{N}(0,1)$. Then,
	under \ref{SAE} with \ref{A}, \ref{C}, \ref{C1} and \ref{C2},
	\begin{align*}
		\frac{1}{m}\sum_{i=1}^m \mbox{P}(|T_i|\leq z_{1-\alpha/2}|\mathbf{v}) = 1-\alpha + O\big(m^{-1/2}\big).
	\end{align*}
\end{prop}

Although nominal coverage is almost surely not attained for single confidence intervals, the coverage probability of marginal confidence intervals under the conditional law still attains its nominal level on average over all clusters, compare also with \citet{Zhang2007}. For the simulated data in Figure \ref{cov_fig} the average coverage is 95.4\% (left) and 94.9\% (right).
The finding in Proposition \ref{avgcov} has been previously described by \citet{Wahba1983} and \citet{Nychka1988} in the context of nonparametric regression.
For an extension of Proposition \ref{avgcov}, more simulation results, and the
construction of Tukey's Intervals, see our Supplement.

\begin{table}[t]
	\centering
	\caption{\label{condtble} \textmd{Coverage of 95\%-confidence ellipsoids in model (\ref{NER}) under conditional law. The relative size to the marginal REML based sets is given in brackets.
	}}
	\vspace*{5pt}
	\begin{tabular}{c rr |cc c |cc }
		&&& Marginal & \multicolumn{2}{c|}{Conditional} &Marginal &Conditional \\
		&$m$&$n_i$ & known $\dlta^v$&known $\lambda$, $\dlta^v$& known $\dlta^v$&REML&REML \\
		\hline\multirow{5}{*}{\shortstack[l]{$\sigma_v^2 = 0.8$\\$ \sigma_e^2 = 0.4$}}
		&5&5&		0.929 \scriptsize{(1)}&0.929 \scriptsize{(1.00)}&		0.928 \scriptsize{(1.01)}&0.895 \scriptsize{(1)}&0.872 \scriptsize{(0.82)}\\
		&50&5&	0.943 \scriptsize{(1)}&0.943 \scriptsize{(0.97)}&		0.946 \scriptsize{(1.26)}&0.921 \scriptsize{(1)}&0.917 \scriptsize{(0.88)}\\
		&5&10&	0.940 \scriptsize{(1)}&0.940 \scriptsize{(1.00)}&		0.940 \scriptsize{(1.00)}&0.922 \scriptsize{(1)}&0.916 \scriptsize{(0.94)}\\
		&50&10&	0.948 \scriptsize{(1)}&0.948 \scriptsize{(1.00)}&		0.948 \scriptsize{(1.03)}&0.938 \scriptsize{(1)}&0.937 \scriptsize{(0.97)}\\
		&10&50&	0.947 \scriptsize{(1)}&0.947 \scriptsize{(1.00)}&		0.947 \scriptsize{(1.00)}&0.944 \scriptsize{(1)}&0.944 \scriptsize{(1.00)}\\
		\hline\multirow{5}{*}{\shortstack[l]{$\sigma_v^2 = 0.6$\\$ \sigma_e^2 = 0.6$}}
		&5&5&		0.926 \scriptsize{(1)}&0.927 \scriptsize{(1.00)}&		0.923 \scriptsize{(1.04)}&0.895 \scriptsize{(1)}&0.840 \scriptsize{(0.67)}\\
		&50&5&	0.944 \scriptsize{(1)}&0.943 \scriptsize{(0.92)}&		0.955 \scriptsize{(2.56)}&0.919 \scriptsize{(1)}&0.924 \scriptsize{(1.13)}\\
		&5&10&	0.938 \scriptsize{(1)}&0.939 \scriptsize{(1.00)}&		0.938 \scriptsize{(1.01)}&0.920 \scriptsize{(1)}&0.906 \scriptsize{(0.87)}\\
		&50&10&	0.948 \scriptsize{(1)}&0.948 \scriptsize{(0.99)}&		0.950 \scriptsize{(1.12)}&0.938 \scriptsize{(1)}&0.938 \scriptsize{(0.98)}\\
		&10&50&	0.947 \scriptsize{(1)}&0.946 \scriptsize{(1.00)}&		0.947 \scriptsize{(1.00)}&0.944 \scriptsize{(1)}&0.944 \scriptsize{(1.00)}\\
		\hline\multirow{5}{*}{\shortstack[l]{$\sigma_v^2 = 0.4$\\$ \sigma_e^2 = 0.8$}}
		&5&5&		0.921 \scriptsize{(1)}&0.929 \scriptsize{(1.02)}&		0.921 \scriptsize{(1.16)}&0.894 \scriptsize{(1)}&0.766 \scriptsize{(0.46)}\\
		&50&5&	0.945 \scriptsize{(1)}&0.942 \scriptsize{(0.79)}&		0.975 \scriptsize{(31.0)}&0.900 \scriptsize{(1)}&0.949 \scriptsize{(14.2)}\\
		&5&10&	0.935 \scriptsize{(1)}&0.942 \scriptsize{(1.01)}&		0.936 \scriptsize{(1.06)}&0.915 \scriptsize{(1)}&0.881 \scriptsize{(0.76)}\\
		&50&10&	0.947 \scriptsize{(1)}&0.947 \scriptsize{(0.96)}&		0.952 \scriptsize{(1.57)}&0.936 \scriptsize{(1)}&0.940 \scriptsize{(1.14)}\\
		&10&50&	0.947 \scriptsize{(1)}&0.947 \scriptsize{(1.00)}&		0.947 \scriptsize{(1.00)}&0.944 \scriptsize{(1)}&0.944 \scriptsize{(1.00)}\\
	\end{tabular}
\end{table}

For the multiple inference we consider the same design as above.
The cluster sample sizes vary from $n_i=5$, cf.\ the study of \cite{Battese1988}, $n_i = 10$ to $n_i = 50$.
A study for an unbalanced data set, where in some cluster even $n_i = 1$ is taken, is provided in the Supplement.
We investigate different ratios of $\sigma_v^2$ and $\sigma_e^2$ with values that again are motivated by the case study estimates in Section \ref{RealStud}. To study different ratios instead of only the one that appears in our case study is interesting
because for model (\ref{NER}) the BLUP can be expressed as an average of a direct estimator with weight $\gamma_i$ and
an estimator for the effect shared by all clusters (often called national estimator) with weight $1-\gamma_i$. Often, $\gamma_i$  is referred to as intraclass correlation coefficient (ICC)
\begin{align*}
	\gamma_i = \frac{\sigma_v^2}{\sigma_v^2 + \sigma_e^2/n_i} \ ,
\end{align*}
and plays a key role in the reliability of $\widehat\mu_i$.
Estimates $\hat{\MU}$ and $\widehat{\SIG}$, as well as $\widehat{\SIG}_v$ and $\hat{\lambda}$ are calculated, and it is checked whether $\MU$ lies within the $95\%-$confidence sets $\Ma$ and $\Ca$.
The structure of $\widehat\SIG$ and $\widehat\SIG_v$ allows matrix inversion by Woodbury's formula, leading to fast calculations.
Table \ref{condtble} contains results based on $10,000$ Monte Carlo samples.
Coverage probabilities are reported together with those of the oracle confidence sets for known $\dlta^v$ and $\lambda$.
Since the simulation is carried out in the conditional setting, $\dlta^v$ is the adequate oracle for the marginal set as well.
By construction, $\hat\dlta$ is a consistent estimator thereof.
The relative average volume of the confidence sets to the volume of the marginal set is given in brackets.
Recall that the asymptotic behavior relies on $m$.

Table \ref{condtble} gives the empirical coverage of the corresponding confidence sets for the nominal level $1-\alpha = 0.95$.
The two most right columns compare the confidence sets based on $\Ma$ and $\Ca$ constructed as outlined above when $\hat\dlta$
is obtained by REML. The other three columns, in the center of the same table, are given for comparison:
The exact confidence set is the `Conditional: known $\lambda$, $\dlta^v$'. Here, nominal level is readily attained.
The coverage of the confidence set `Marginal: known $\dlta^v$' exhibits the error solely due to the
misspecification in using the marginal set for the conditional setting, as described in Theorem \ref{margforcond}. Despite its error rate $O(m^{-1/2})$, the empirical coverage is so close to $1-\alpha$, that it cannot be seen in Table \ref{condtble}.
Comparing the two `Conditional' columns with oracle parameters reveals the impact of estimating $\lambda$; see the Supplement for a deeper analysis of the reliability of the estimation of $\lambda$.
%
\begin{figure}[tb]
	    \vspace{25pt}
	    	\centering
		\input{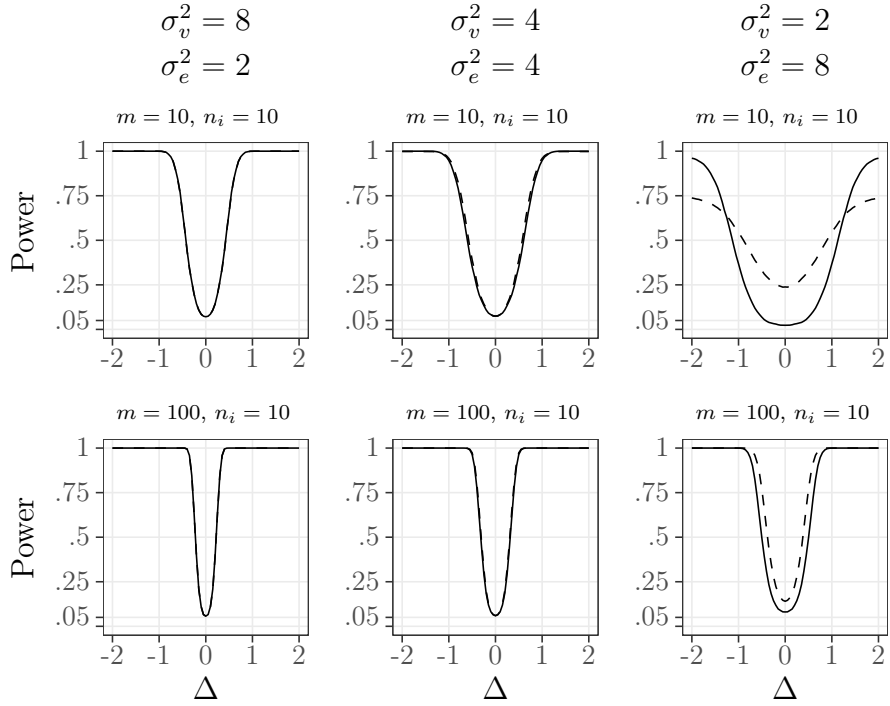}
		\vspace*{5pt}
		\caption{Power of tests based on confidence ellipsoids $\Ma$ (solid line) and $\Ca$ (dashed) for model (\ref{NER}) in the conditional setting with $H_1:  \MU = \mathbf{a} + \mathbf{1}_m\Delta$.}
		\label{power}
\end{figure}
Clearly, the coverage probabilities improve for larger $m$ and/or $n_i$.
This is in line with the theoretical findings.
However, the coverage error is superimposed by the shape of the ICC.
If it is close to $1$, the REML estimates are stable, and similarly $\hat\lambda$.
The ICC is influenced by two drivers:
Firstly, by the relative size of $\sigma_v^2$ to $\sigma_e^2$.
If $\sigma_v^2$ is large, the empirical coverage is closest to the nominal level.
This has already been observed for individual confidence intervals \citep{Das2004}.
Secondly, the ICC relies on the size of $n_i$.
Irrespective of the reliability of the REML estimates, a large $n_i$ results in accurate coverage probabilities, as can be seen on the last row for each configuration of $(\sigma_v^2, \sigma_e^2)$.
Conversely, even for known $\dlta^v$, a small $n_i$ may cause in a severe under-coverage.
All these effects shape the performance of the REML based confidence sets:
it is evident that the asymptotic behavior cannot be observed when $n_i$ and $\sigma_v^2$ are small compared to the noise level $\sigma_e^2$. For more discussion see the Supplement.


Finally we consider the test $H_0: \MU = \mathbf{a}$ vs.\ $H_1:  \MU = \mathbf{a} + \mathbf{1}_m\Delta$, $\mathbf{a}\in\mathbb{R}^m$ with $\Delta\in\mathbb{R}$. Power functions studying the error of the second kind for different parameters $m$ and $n_i$ are given in Figure \ref{power} for different ICC.
%
Unsurprisingly, the power growths steeper for larger $m$ and $n_i$, but again is sensitive to the relative size of $\sigma_v^2$ to $\sigma_e^2$. The power of the tests based on marginal sets (solid line) is notably steeper than the slope of the power based on conditional sets (dashed).

All in all, both the conditional and marginal sets exhibit similar coverage probabilities, which is in line with the theoretical findings.
However, due to its simpler construction and broader application, the results of both Table \ref{condtble} and Figure \ref{power} favor the use of marginal confidence sets, especially for testing.

\section{Study on Covid-19 Mortality in US State Prisons}\label{RealStud}

The methods introduced above are applied to Covid-19 related mortality rates in US state prisons between
March 2020 until the end of March 2021, published by \citeauthor{Covid19}.
The data are from $n = 494$ US state prisons of $m = 45$ states which form the clusters.
The model is
$ y_{ij} = \beta_0 + x_{ij}\beta_1 + v_i + e_{ij}  ,
$
where $y_{ij}$ is the log-mortality for prison $j$ in state $i$ and $x_{ij}$ the standardized county log-mortality in which the prison is located.
The covariates account for local effects on mortality, while the
error terms account for the plethora of unobserved variables.
The random effect $v_i$ describes the remaining state effect on mortality.
The number of prisons $n_i$ in each state ranges from $1$ to $46$, with a median of $8$.
The use of direct estimators is unreliable due to the small number of observations per state,
so that the LMM modeling device is appealing. However, following the logic of our above discussions, the inference will be performed conditionally on $v_i$. The parameter of interest $\mu_i = \beta_0 + \sum_{j=1}^{n_i} x_{ij} \beta_1 / n_i + v_i$, the mean log-mortality in prisons per state (subsequently ``mortality''), is estimated via the EBLUP $\widehat\mu_i$. The fixed effects are estimated as $(\widehat\beta_0, \widehat\beta_1)\approx (-4.79, 0.03)$ and the variance components via REML as $\widehat\sigma_v^2 \approx 0.43$ and $\widehat\sigma_e^2 \approx 0.86$. In sum, the setting is similar to those analyzed in our simulation studies. Furthermore, the specification, including the normality assumption for model (\ref{LMM}), is graphically assessed by residual analyses in the Supplement. The residuals show neither any anomalies, nor a violation of our model assumptions.
The estimates $\widehat\mu_i$ are visualized in Figure \ref{usa}.

\begin{figure}
	    	\centering
				\begingroup
				\input{usa2.tex}
				\endgroup
		\vspace*{-25pt}
		\caption{Conditional state means $\widehat\mu_i$ of Covid-19 mortality in US state prisons.}
		\label{usa}
\end{figure}

First, we state the hypothesis that the state effect is due to Covid-19 related policies.
An interesting question could be if mortality in democratic governed states is lower than in republican ones.
Formally, let $\mu_R$ be the mortality for all $22$ states governed by republicans and $\mu_D$ for all $23$ democratic ones, for which data is available.
The corresponding t-test to $H_0$: $\mu_R \leq \mu_D$ vs. $H_1$: $\mu_R > \mu_D$ using direct estimates for the two types of states, rejects the null for common significance levels with a p-value of $P_{H_0}(T>t) \approx 10^{-8} $.
However, the above t-test supposes that observations given the same party come all from a distribution with the same mean, i.e., there were no systematic differences in mortality within democratic or republican states, respectively.
This can be checked using a linear hypothesis test as described in Section \ref{SimInf}.
Formally, interest lies in verifying the hypothesis that groups of states share the same state effect.
For a group of $u + 1$ states, let $\mathbf{L} = (\boldsymbol{0}_{u}, \dots,  \mathbf{L}^{\ast}, \boldsymbol{0}_{u}, \dots)$ be  $(u \times 45)$, with $\mathbf{L}^{\ast} = \big(\mathbf{I}_{u}, \boldsymbol{0}_{u}\big) - \mathbf{1}_{u}\mathbf{1}_{u + 1}^t/ (u + 1)$ corresponding to the states of interest.
We test the hypothesis $H_0: \mathbf{L}\MU = \boldsymbol{0}_{u}$ against $H_1: \mathbf{L}\MU \neq \boldsymbol{0}_{u}$.
As described in Section \ref{SimInf}, this tests whether all states in the considered group share an equal state mean.
If the group consists of the first $u + 1$ states, the null hypothesis is equivalent to $H_0: \mu_1 = \mu_2 = \dots = \mu_{u+1}$.
The result for two such tests for equality for democratic and republican governed states respectively, are given in Table \ref{governor}.
The table reports the rank of $\mathbf{L}$, i.e., the rate at which the tests are consistent, the value of the quadratic form as `Pivot'
together with the corresponding quantile and p-value and -- for the conditional case -- the estimated non-centrality parameter $\hat{\lambda}_{\mathbf{L}}$.

\begin{table}
\centering
\caption{\label{governor} \textmd{Tests for the equality of state means of groups by governor.
}}
\vspace{5pt}
\begin{tabular}{l c r c  l r c r l}
\hline
\multirow{2}{*}{Governor}&\multirow{2}{*}{$u+1$}& \multicolumn{3}{c}{Marginal} & \multicolumn{4}{c}{Conditional} \\
&& \small{Pivot}&
\small{$\chi^2_{u+1, 0.95}$}&
\small{p-value}&
\small{Pivot}&
\small{$\chi^2_{u+1, 0.95}(\hat\lambda_{\mathbf{L}})$}&
\small{$\hat\lambda_{\mathbf{L}}$}&
\small{p-value}\\
\hline
Democrat&$23$&$77$ &$34$ & $5\times10^{-8}$& $104$ &$87$ &$40$ & $5\times10^{-3}$\\
Republican&$22$&$141$ &$33$ & $8\times10^{-20}$& $172$ &$35$ &$1$ & $6\times10^{-24}$\\
\hline
\end{tabular}
\end{table}

At common significance levels, both tests reject the hypothesis that the mortality is equal in all democratic or republican states, respectively, which conflicts the assumptions of the above t-test.
It is however noteworthy, that the $p$-values of the marginal and conditional set are noticeably different.
This is due to the estimate obtained for the non-centrality parameter, on which the conditional confidence set relies.
The estimator for the non-centraility parameter $\lambda$ works well in balanced panels, as can be seen in Table \ref{condtble}.
In the present data set, however, the panel is unbalanced with the large proportion of states having few observations. Moreover, the large error variance $\sigma_e^2$ in relation to $\sigma_v^2$ in such unbalanced panels with small sample sizes is known to lead to very unreliable estimators for $\sigma_e^2$ and $\sigma_v^2$ and herewith for $\lambda_{\mathbf{L}}$, especially if $m$ is not too large.
This is confirmed by parametric bootstrap estimates for $\lambda_{\mathbf{L}}$ in all  tests of Tables \ref{governor} and \ref{census}, provided in the Supplement.
The variable estimator $\lambda_{\mathbf{L}}$ causes the conditional confidence set
to be unreasonably large, resulting in overcoverage. This claim is verified in the last line of Table 6 in the Supplement, which replicates exactly the setting corresponding to our real data example.  At the same time, the results shown in the last line of Table 6 confirm that the marginal sets perform excellent in the given parameter constellation.
Moreover, if one is interested in other than the previously considered groups, the conditional approach requires to re-estimate the non-centrality parameter on each new subset of interest.
These aspects enhance the value of Theorem \ref{margforcond} and give a strong support for application of the marginal set in practice.

Instead of looking at political party effects, one may look at geographic effects, and check if among certain groups of states their mortality is equal. We repeat the above test for groups formed by the four regions of the US census bureau.
The results are given in Table \ref{census}.
For common significance levels, the tests reject the null hypothesis for the census regions Midwest, South and West.
For the census region Northeast, the null hypothesis cannot be rejected.
Potentially, this is because the state policies are homogeneous within this census region.
Again, the influence of the non-centrality parameter can be observed for Northeast, even though it may not make a difference for the conclusion as
the marginal and conditional tests give the same results for significance levels $\alpha = 0.01$, $0.05$, and $0.1$.
For the southern census region, $10$ of $15$ of all individual tests do not reject $H_0$, and neither would a joint test with Bonferroni correction.
In fact, the latter is true for all census regions except Northeast.
This illustrates that our multiple test represent an important complement to the existing single ones, combined or not with Bonferroni.

\begin{table}
\centering
\caption{\label{census} \textmd{Tests for the equality of state means of groups by census regions.
}}
\vspace{5pt}
\begin{tabular}{l c r c  l r c r l}
\hline \vspace*{2pt}
\multirow{2}{*}{\shortstack{Census\\region}}&\multirow{2}{*}{$u+1$}& \multicolumn{3}{c}{Marginal} & \multicolumn{4}{c}{Conditional} \\
&& \small{Pivot}&
\small{$\chi^2_{u+1, 0.95}$}&
\small{p-value}&
\small{Pivot}&
\small{$\chi^2_{u+1, 0.95}(\hat\lambda_{\mathbf{L}})$}&
\small{$\hat\lambda_{\mathbf{L}}$}&
\small{p-value}\\
\hline
\hline
Midwest&$10$ &$45$ &$17$ &$8\times10^{-7}$&$53$ &$17$&$0.3$ & $5\times10^{-8}$\\
Northeast&$8$&$8$ &$14$&$0.33$&$19$ &$71$ &$41$ & $0.99$\\
South&$16$&$131$ &$25$&$1\times10^{-20}$&$147$ &$25$&$0$ & $8\times10^{-24}$\\
West&$10$&$19$ &$17$&$2\times10^{-2}$&$23$ &$17$&$0$ & $5\times10^{-3}$\\
\hline
\end{tabular}
\end{table}

Certainly, the above illustration gives just some particular examples, but it is obvious that
any other linear hypotheses with $u\leq m$ could be tested analogously.
We believe that such tests are highly relevant, insightful and helpful in practice.
One can also use the confidence sets to see in which clusters one needs to change how much in order to eliminate significant differences, employing thereby our tools for policy makers.


\section{Discussion} \label{Discussion}

Under assumption \ref{SAE} inference based on predictors for single clusters is intractable under conditional law due to the bias.
This is the reason why single cluster inference has only been performed under the marginal framework.
As shown in Proposition \ref{avgcov}, the inference for the individual mixed parameter holds on average only.
In this work 
we derived joint confidence sets for mixed parameters $\mu_1,\ldots,\mu_m$ in LMMs under both, marginal and
conditional law. The latter requires the estimation of a non-centrality parameter of the respective $\chi^2(\lambda)$-distribution.
We have shown that with its estimate, the desired nominal coverage is attained at the usual parametric rate.
To the best of our knowledge, our method allows for inference on multiple clusters under the conditional law for the first time.
In particular, it lends itself to infer on a subset of clusters of interest, as illustrated in the study on the Covid-19 mortality in US state prisons.
Further, we show that, surprisingly and in contrast to cluster-wise confidence intervals, the joint (or multiple) confidence sets built under marginal law are approximately valid at the same parametric rate when conditioning on the clusters.
A simulation study confirms this effect already for samples of small and moderate size.
Our results hold for all kind of linear combinations of mixed parameters $\mu_i$ of a cluster $i$.

The order of the derived error relies on the normality assumption in (\ref{LMM}).
If no distributional assumption is justified, additional regularity conditions governing the boundedness of higher moments have to be imposed, and resampling methods could be applied.
Moreover, simulations carried out for non-Gaussian random effects, shown in the Supplement, indicate the robustness of the proposed confidence sets.
Furthermore, when it is of interest to test linear contrasts of mixed parameters, we extend our test to cover multiple comparisons by Tukey's method, see the Supplement.
However, the application of this method is limited to special cases where the corresponding bias can be shown to be negligible, and the considered subset of pairwise differences falls exactly into the class of Tukey's testing problem.
Finally, we expect that generally, our methods and results can be extended
to other predictors of LMMs, such as the empirical best predictor of \citet{Jiang2011}.

%
%


\begin{appendix}
\label{appendix}

 \section*{Acknowledgments}
{The authors thank Domingo Morales, Carmen Cadarso-Su\'arez, Jiming Jiang, Mar\'{\i}a-Jos\'e Lombard\'{\i}a and Wenceslao Gonz\'alez-Manteiga for helpful discussion. 
This work has been carried out while the first two authors were employed at the University of Göttingen, Germany. They also acknowledge the funding by the German Research Association (DFG) via Research Training Group 1644 ``Scaling Problems in Statistics''; the last author acknowledges financial support from the Swiss National Science Foundation, project 200021-192345.}

\section*{Appendix}

\subsection*{Asymptotic scenarios beyond \ref{SAE}}

Although the results are derived under \ref{SAE}, they are not restricted to such an asymptotic scenario.
Under the marginal law, $\mbox{E}(\hat\MU - \MU) = \mathbf{0}_m$, but under \ref{SAE}, $\mbox{E}(\hat\MU - \MU|\mathbf{v}) \nrightarrow \mathbf{0}_m$ under the conditional law.
If \ref{SAE} were to be relaxed, and both $m\rightarrow\infty$ and all $n_i\rightarrow\infty$, the EBLUP would be consistent under both probability measures,
and both, marginal and conditional, collapse into one.
To investigate the effect of unbounded $n_i$ on Theorems \ref{margdist}-\ref{margforcond}, the source of error terms becomes crucial.
The error term in Theorem \ref{margdist} is due to $\hat\dlta = \dlta + \{O(m^{-1/2})\}_r$ while the one in Theorem \ref{conddist} is due to $\hat\dlta = \dlta^v + \{O(m^{-1/2})\}_r$ and the estimation of $\lambda$.
As the EBLUP is consistent under conditional law if $n_i\rightarrow\infty$ for all $i$, the non-centrality parameter vanishes in such cases, i.e., $\lambda \rightarrow 0$.
The same holds for Theorem \ref{margforcond}.
%
Technically, the cases for unbounded $n_i$ differ from \ref{SAE} as
the leading entries on the diagonal of $\SIG$ and $\SIG_v$ vanish: $(\SIG)_{ii} = O(n_i^{-1})$ and $(\SIG_v)_{ii} = O(n_i^{-1})$.
In order to assess which asymptotic behavior (of the diagonal entries or $\hat\dlta$) determines the rate of the error term in each theorem, it is required to fix the relation of $m$ and $n_i$, $i=1,\dots, m$.

Under the asymptotic scenario $m\rightarrow\infty$ and all or some $n_i\rightarrow\infty$,
the stated results still hold, and the error rates can improve.
This depends intricately on the number of unbounded cluster sample sizes and the rate at which they grow.
If the number of clusters with bounded sample sizes is itself unbounded, that is $O(m)$, the error rates generally fall back on what is stated in Theorems \ref{margdist}-\ref{margforcond}.
If it is bounded, a toy example shows that they can improve.
Set the sample size of a single cluster as fixed and let all other cluster sample sizes grow at the same rate as the number of clusters $m$.
That is, $n_1 = O(1)$, $m\rightarrow\infty$ and $m/n_i \rightarrow 1$, for $i = 2, \dots, m$.
Then, by Lemma \ref{marglemma} and the proof of Theorem \ref{margdist}, $m^{-1/2}\sum_{i=1}^m(\hat\mu_i - \mu_i)^2 = O_p(m^{-1/2})$, so that the error in Theorem \ref{margdist} is reduced to $O(m^{-3/2})$.

\subsection*{Proofs}

The proofs are given in two parts.
First, the order of the bias of the covariance matrix estimator is established in Lemmas \ref{marglemma} and \ref{LemmaCondReml}.
The former is omitted for brevity, the latter given in the Supplement.
Both rely on Taylor approximations, similar to \citet{Prasad1990} and \citet{Datta2000}.
The difficulty in Lemma \ref{LemmaCondReml} lays in the nature of $\dlta^v$ and decomposition (\ref{SIGCorig}), a multitude of additional terms have to be evaluated.
In the second part of proofs it is shown that the resulting error rate is preserved
in the evaluations that lead to Theorems \ref{margdist}-\ref{margforcond}.
Since both the dimension of the covariance matrix estimator as well as the error rate are given in terms of $m$, this has to be carefully addressed in matrix inversion, the quadratic form, and the final probabilistic statement.

\subsubsection*{Proof for Theorem \ref{margdist}}
\begin{proof} 
Let $(\SIG)_{ik} = \sigma_{ik}$, $\{\widehat\SIG(\hat\dlta)\}_{ik} = \hat\sigma_{ik}$ and $\{\widehat\SIG(\dlta)\}_{ik} = \tilde\sigma_{ik}$.
We first show that
\begin{align}\label{firsttoshow}
\|\widehat{\SIG}^{-1/2}(\hat{\MU} - \MU)\|^2 &=  \|\SIG^{-1/2}(\hat{\MU} - \MU)\|^2 + O_p(m^{1/2}).
\end{align}
By Lemma \ref{marglemma}, $\mbox{E}(\hat\sigma_{ik}) = \sigma_{ik} + O_p(m^{-3/2})$, as well as $\tilde\sigma_{ik}= \sigma_{ik} + O_p(m^{-3/2})$. Note that $\hat\delta_e - \delta_e = O_p(m^{-1/2})$. Further, $\tilde\sigma_{ii}=O(1)$ as well as $\tilde\sigma_{ik}=O(m^{-1})$ for $i\neq k$ and this order is preserved for its derivatives with respect to $\dlta$. Thus,
\begin{align*}
\mbox{Var}(\hat{\sigma}_{ik}) &= \mbox{E}\big[ \{\hat{\sigma}_{ik} - \tilde{\sigma}_{ik} + O(m^{-3/2}) \}^2\big]  \\
&= \mbox{E}\bigg[\bigg\{ (\hat\dlta - \dlta)^t\frac{\partial\tilde{\sigma}_{ik}}{\partial\dlta} + (\hat\dlta - \dlta)^t\frac{\partial^2\tilde{\sigma}_{ik}}{\partial\dlta\partial\dlta^t}(\hat\dlta - \dlta) + O_p(m^{-3/2}) \bigg\}^2\bigg] \\
&= \mathbbm{I}(i=k)O(m^{-1}) + O(m^{-3}).
\end{align*}
By Chebychevs inequality, for a random variable $X$ with finite variance $X = \mbox{E}(X) + O_p\{\sqrt{\mbox{Var}(X)}\}$.
It follows that $\widehat\SIG = \SIG - \C$ where
\begin{align*}
\SIG &= \mbox{diag}[\{O(1)\}_m] + \{O(m^{-1})\}_{m \times m}, \\
\C &= \mbox{diag}[\{O_p(m^{-1/2})\}_m] + \{O_p(m^{-3/2})\}_{m \times m}.
\end{align*}
It is now shown that inverting preserves the error.
Let $\D = (\dd_1, \dots, \dd_m)$ for $\dd_i$ as in \ref{C} and note that $(\X^t\V^{-1}\X)^{-1} = \{O_p(m^{-1})\}_{p\times p}$.
The matrix inversion formula yields
\begin{align*}
\SIG^{-1} 	&= \big\{  \mathbf{K}_1 + \D^t\big(\X^t\V^{-1}\X\big)^{-1}\D   \big\}^{-1}  \\
			&= \mathbf{K}_1^{-1} - \mathbf{K}_1^{-1}\D^t\big(    \X^t\V^{-1}\X       + \D\mathbf{K}_1^{-1}\D^t \big)^{-1}\D\mathbf{K}_1^{-1}  \
			= \mathbf{K}_1^{-1} + \{O(m^{-1})\}_{m\times m}.
\end{align*}
Thus, $\C\SIG^{-1} = \mbox{diag}[\{O_p(m^{-1/2})\}_m] + \{O_p(m^{-3/2})\}_{m \times m}$. Denote $\lambda_{\C\SIG^{-1}}$ as largest eigenvalue of $\C\SIG^{-1}$. With the column-sum norm, $\lambda_{\C\SIG^{-1}} \leq \max_{k=1,\dots,m} \sum_{i=1}^m | \{\C\SIG^{-1}\}_{ik} | = O(m^{-1/2}) < 1$ for large $m$. Writing the inverse as Neumann-series, $( \I_m  - \C\SIG^{-1}    )^{-1} =  \I_m + \mbox{diag}[\{O_p(m^{-1/2})\}_m] + \{O_p(m^{-3/2})\}_{m \times m}$. Now
\begin{align*}
\widehat\SIG^{-1} 	&= \SIG^{-1}\big( \I_m - \C\SIG^{-1}\big)^{-1} = \SIG^{-1} +  \mbox{diag}[\{O_p(m^{-1/2})\}_m] + \{O_p(m^{-3/2})\}_{m \times m}
\end{align*}
Eventually, since $m^{-1/2}\sum_{i=1}^m(\hat\mu_i - \mu_i)^2 = O_p(m^{1/2})$ and
$Q = \|\SIG^{-1/2}(\hat{\MU} - \MU)\|^2/m = O_p(1)$, putting all parts together gives (\ref{firsttoshow}).
Further, let $U=\|\widehat{\SIG}^{-1/2}(\hat{\MU} - \MU)\|^2/m -  \|\SIG^{-1/2}(\hat{\MU} - \MU)\|^2/m=O_p(m^{-1/2})$ with probability density function $f_U$ and let $z=m^{-1}\chi^2_{m, 1-\alpha} =O(1)$, such that
\begin{eqnarray*}
&&\mbox{P}\bigg\{ \|\widehat{\SIG}^{-1/2}(\hat{\MU} - \MU)\|^2 < \chi^2_{m, 1-\alpha} \bigg\} = \mbox{P}\big( Q + U < z\big)
    =  \\
&& \int_{\mathbb{R}}  \mbox{P}\big( Q < z - u\big) f_U(u)  du
= \int_{\mathbb{R}}  \bigg\{ \mbox{P}\big( Q < z \big) + O(m^{-1/2}) \bigg\} f_U(u)  du   
= 1-\alpha + O(m^{-1/2})  ,
\end{eqnarray*}
which concludes the proof.
\end{proof}

\subsubsection*{Proof and Definitions for Theorem \ref{conddist}}


First, define $\mathbf{w}_i = (\mathbf{b}_i^t\mathbf{Z}_i - \mathbf{h}_i^t)\mathbf{J}_i + \mathbf{d}_i^t(\X^t\V^{-1}\X)^{-1}\X^t\V^{-1} \in \mathbb{R}^n$, $n = \sum_{i=1}^m n_i$, so that $\mathbf{w}_i^t\mathbf{e} = \tilde{\mu}_i - \mbox{E}(\tilde{\mu}_i|\mathbf{v})$.
Let ${\mathbf{L}}_3^{\ast}(\dlta^v) = \mbox{Cov}\left(
\hat\MU-\tilde\MU, \tilde\MU
\big| \mathbf{v}\right)$ and $\widehat{\mathbf{L}}_3(\dlta^v) = \widehat{\mathbf{L}}_3^{\ast}(\dlta^v)
 + \widehat{\mathbf{L}}_3^{\ast}(\dlta^v)^t$.
If $\dlta^v$ is estimated via
\begin{enumerate}[label=(\roman*),font=\upshape, itemsep=1ex]
\item REML, given $\mathbf{P} = \V^{-1} -\V^{-1}\X(\X^t\V^{-1}\X)^{-1}\X^t\V^{-1}$, then
\begin{equation}
\begin{aligned}\label{L3REML}
\widehat{\mathbf{L}}_3^{\ast}(\dlta^v) &= \left[\rule{0cm}{.8cm}\right.
2\sum_{e=1}^r\text{tr}\bigg\{\mathbf{P}
\frac{\partial\mathbf{V}}{\partial\delta^v_e}\mathbf{P}\R   \mathbf{w}_i( \overline{\V})_e^t\frac{\partial\mathbf{w}_k^t}{\partial\dlta^v}  \R \bigg\}  \\
&  \hspace*{2.5cm} + 4\displaystyle\sum_{e,d=1}^r \text{tr}\left\{
\sum_{f,g=1}^r ( \overline{\V})_{ef}( \overline{\V})_{fg} \mathbf{w}_i\frac{\partial^2\mathbf{w}_k^t}{\partial\delta^v_e\partial\delta^v_d}
\R \right\}\big(\overline{\V}^{-1}\big)_{ed}\\
&  \hspace*{2.5cm}-2\displaystyle\sum_{e,d=1}^r \text{tr}\left\{
\sum_{f=1}^r ( \overline{\V})_{ef}\mathbf{w}_i( \overline{\V})_d^t
\frac{\partial\overline{\V}^{-1}}{\partial\delta^v_e} \overline{\V}
\frac{\partial\mathbf{w}_k^t}{\partial\dlta^v}
\R \right\}\big(\overline{\V}^{-1}\big)_{ed}\\
&\hspace*{2.5cm} + 2\displaystyle\sum_{e,d,g=1}^r \text{tr}\bigg\{\mathbf{w}_i( \overline{\V})_e^t\frac{\partial\mathbf{w}_k^t}{\partial\dlta^v}\R \bigg\}\frac{\partial(\overline{\V}^{-1})_{ed}}{\partial\delta^v_g}(\overline{\V})_{ed}
\left]\rule{0cm}{.8cm}\right._{i,k=1\dots, m}.
\end{aligned}
\end{equation}
\item Henderson III, then
\begin{align} \label{L3H3}
\widehat{\mathbf{L}}_3^{\ast}(\dlta^v) =
\left[ \sum_{e=1}^r
  2\text{tr}\left\{\mathbf{w}_i\frac{\partial\mathbf{w}_k^t}{\partial\delta^v_e} \R\mathbf{C}_e\R\right\}
  + \sum_{g=1}^r \text{tr}\left\{\mathbf{w}_i\frac{\partial^2\mathbf{w}_k^t}{\partial\delta^v_e\partial\delta^v_g}\R\right\} \overline{\V}_{eg} \right]_{i,k=1\dots,m}.
\end{align}
\end{enumerate}
Both estimators have entries of order $O(m^{-1})$.
Their are derived in analogy to $\widehat{\mathbf{L}}_4(\dlta^v)$, which is outlined in the Supplement.


\begin{proof}
With Lemma \ref{LemmaCondReml} the proof for Theorem \ref{margdist} can be replicated giving
\begin{align*}
\mbox{P}\bigg\{ \|\widehat{\SIG}_v^{-1/2}(\hat{\MU} - \MU)\|^2 < \chi^2_{m, 1-\alpha}(\lambda) \bigg|\vv \bigg\} &= 1-\alpha + O(m^{-1/2}).
\end{align*}
It thus remains to show that
\begin{align*}
\chi^2_{m, 1-\alpha}(\lambda) = \chi^2_{m, 1-\alpha}\left\{\tilde\lambda\left( \widehat\SIG_v, \hat\bta, \hat\dlta  \right)\right\} + O_p\left(m^{1/2}\right).
\end{align*}
Examining the entries of both $\mathbf{A}$ and $\SIG_v^{-1}$ gives, using the decomposition of the proof of Theorem \ref{margdist},
$\widehat\SIG_v^{-1} = \SIG_v^{-1} + \mathbf{B}$, for
$\mathbf{B} = \text{diag}[\{O(m^{-1/2})\}_{n_i \times n_i}] + \{O(m^{-3/2})\}_{n\times n}$, that
 $\mathbf{A}^t\SIG_v^{-1}\mathbf{A} = \text{diag}[\{O(1)\}_{n_i \times n_i}] + \{O(m^{-1})\}_{n\times n}$, so that
 \begin{align*}
 	\mathbf{A}^t\widehat\SIG_v^{-1}\mathbf{A}
	= \mathbf{A}^t\SIG_v^{-1}\mathbf{A}
	+ \text{diag}[\{O(m^{-1/2})\}_{n_i \times n_i}] + \{O(m^{-3/2})\}_{n\times n}.
 \end{align*}
Using \ref{C1} and \ref{C2}, so that $\hat{\dlta} = \dlta^v + \{O_p(m^{-1/2})\}_r$ as given by Lemma 4 in the Supplement and putting all parts together gives
\begin{align*}
	\tilde\lambda\left( \widehat\SIG_v, \hat\bta, \hat\dlta  \right) = \tilde\lambda\left( \SIG_v, \tilde\bta, \dlta^v  \right) + O_p(m^{1/2}).
\end{align*}
This error rate is sufficient as the estimator effectively contributes as $\hat\lambda/m$ in $\chi^2_m(\hat\lambda)$.
Now we show that $\tilde\lambda = \tilde\lambda( \SIG_v, \tilde\bta, \dlta^v ) = \lambda + O_p(m^{1/2})$ by considering its expectation and variance.
\begin{align*}
	\mbox{E}\big( \tilde\lambda \big| \vv \big)
	 &= \lambda + \|\SIG_v^{-1/2}\mathbf{A}\X(\X^t\V^{-1}\X)^{-1}\X^t\V^{-1}\Z\vv\|^2
	 + 2( \mathbf{A}\X\tilde\bta)^t\SIG_v^{-1}\mathbf{A}\Z\vv\\
	 &= \lambda + O(m^{1/2}),
\end{align*}
using \ref{C1} and \ref{C2}. Similarly, $\text{Var}(\tilde\lambda|\vv)=O(m)$.
Hence, $\tilde{\lambda} = \lambda + O_p(m^{1/2})$. Eventually,
\begin{align*}
\chi^2_{m, 1-\alpha}\big(\hat\lambda\big)
&= \chi^2_{m, 1-\alpha}\big\{\lambda + O_p(m^{1/2})\big\} 
= \chi^2_{m, 1-\alpha}(\lambda)  + O_p(m^{1/2}).
\end{align*}
This concludes the proof.
\end{proof}

\subsubsection*{Proof for Theorem 3}
Another way to obtain a pivotal for multiple inference is to evaluate the distribution of the quadratic form $Q = \|\SIG^{-1/2}(\hat{\MU} - \MU)\|^2$ under the conditional law. It is distributed as generalized non-central $\chi^2$, and thus has no analytically tractable probability density function. However, due to the linearity of $\hat{\MU} - \MU$ in $\mathbf{v}$, the quadratic form $Q$ can be suitably split up in treatable terms.

\begin{proof}
Due to linearity of $\hat{\MU} - \MU$, it holds that 
$\SIG = \SIG_v + \SIG_b$, where $\SIG_b = \mbox{Cov}(\MU_b)$ for $\MU_b = \mbox{E}(\hat{\MU} - \MU|\mathbf{v})$ by the law of total variance. Moreover,
\begin{align*}
\SIG^{-1} &= \big(\SIG_v+\SIG_b\big)^{-1} = \SIG_v^{-1}-\SIG_v^{-1}\big(\SIG_v^{-1}+\SIG_b^{-1}\big)^{-1}\SIG_v^{-1} = \SIG_v^{-1} - \mathbf{T}_c^{-1},
\end{align*}
where $\mathbf{T}_c^{-1}$ fulfills $\SIG_v\mathbf{T}_c^{-1} = \SIG_b\SIG^{-1}$. Now consider $Q= S + R$ with
\begin{align*}
S &= \|\SIG_v^{-1/2}(\hat{\MU} - \MU - \MU_b)\|^2, \\
R &= \|\SIG^{-1/2}\MU_b\|^2 + 2\MU_b^t\SIG^{-1}(\hat{\MU}-\MU-\MU_b) - \|\mathbf{T}_c^{-1/2}(\hat{\MU} - \MU - \MU_b)\|^2.
\end{align*}
It holds that $S|\mathbf{v}\sim\chi^2_m$.
Next, we show that $R$ is of lower order compared to $S$.
Let $\mathbf{W}\in\mathbb{R}^{m\times mq}$ such that $\MU_b=\mathbf{W}\vv$.
Note that $\mathbf{W}= \text{diag}[\{O(1)\}_{1\times q}]_m + \{O(m^{-1})\}_{m\times mq}$.
\begin{align*}
\mbox{E}(R|\vv)
&= \text{tr}\left\{ \SIG^{-1}(\MU_b\MU_b^t - \SIG_b)\right\}
= \text{tr}\left[ \SIG^{-1}\mathbf{W}\{\vv\vv^t-\text{diag}(\G)_m\}\mathbf{W}^t\right] = O(m^{1/2}),
\end{align*}
by the same reasoning as in the proof of Lemma 4 in the Supplement.
Similarly, $\text{Var}(R|\vv) = O(m)$.
Hence, $R = \mbox{E}(R|\vv) + O_p\{\sqrt{\mbox{Var}(R|\vv)}\} = O_p(m^{1/2})$.
Now, using that $S = O_p(m)$,
\begin{align*}
\mbox{P}\big( Q < \chi^2_{m, 1-\alpha} \big|\mathbf{v}\big) &= \mbox{P}\bigg\{ \frac{S}{m} + O_p(m^{-1/2}) < \frac{\chi^2_{m, 1-\alpha}}{m} \bigg|\mathbf{v}\bigg\}\\
&= \mbox{P}\big( S < \chi^2_{m, 1-\alpha} \big|\mathbf{v}\big) + O(m^{-1/2})
= 1-\alpha + O(m^{-1/2}) .
\end{align*}
Replacing $\SIG$ in $Q$ by $\widehat{\SIG} = \SIG + \{O_p(m^{-1/2})\}_{m\times m}$ gives $\|\widehat{\SIG}^{-1/2}(\hat{\MU} - \MU)\|^2/m = Q/m + O_p(m^{-1/2})$ as in the proof of Theorem \ref{margdist}.
The order of the error coincides with the one in above equation, which gives
$\mbox{P}( \|\widehat{\SIG}^{-1/2}(\hat{\MU} - \MU)\|^2 < \chi^2_{m, 1-\alpha} |\mathbf{v})
= 1-\alpha + O(m^{-1/2}) $.
\end{proof}

\end{appendix}

\bibliographystyle{apalike}
\bibliography{References}

\end{document}